\documentclass[11pt]{article}
\usepackage{amsmath}
\usepackage{amssymb}
\usepackage{latexsym}
\usepackage{amsthm}
\usepackage{fullpage}
\usepackage[colorinlistoftodos]{todonotes}
\usepackage{marginnote}
\usepackage{enumitem}
\usepackage{mathtools}
\usepackage{dsfont}
\usepackage{mdframed}
\usepackage{xcolor}
\usepackage{hyperref}
\usepackage{esint}
\usepackage{amsfonts,latexsym,amsmath,amscd,geometry}
\numberwithin{equation}{section}

\newtheorem{theorem}{Theorem}[section]
\newtheorem{lemma}[theorem]{Lemma}

\newtheorem{proposition}[theorem]{Proposition}
\newtheorem{remark}[theorem]{Remark}

\newtheorem{corollary}[theorem]{Corollary}

\definecolor{vine}{rgb}{0.7,0.1,0.1}
\definecolor{vine}{rgb}{0.7,0.1,0.1}
\definecolor{darkgreen}{rgb}{0,0.5,0}

\newcommand{\R}{{\mathbb R}}

\newcommand \nc{\newcommand}

\nc{\ba}{\begin{array}}\nc{\ea}{\end{array}}
\nc{\be}{\begin{eqnarray}}\nc{\ee}{\end{eqnarray}}
\nc{\beq}{\begin{equation}}\nc{\eeq}{\end{equation}}
\nc{\bex}{\begin{eqnarray*}}\nc{\eex}{\end{eqnarray*}}
\nc{\btm}{\begin{theorem}} \nc{\etm}{\end{theorem}}
\nc{\blm}{\begin{lemma}} \nc{\elm}{\end{lemma}}
\nc{\ld}{\lambda}
\nc{\va}{\varphi}
\nc{\ve}{\varepsilon}

\def\pf{\noindent{\bf Proof.\quad}}

\author{Zhiyuan Geng\footnote{Department of Mathematics, Purdue University, West Lafayette, IN 47907 (geng42@purdue.edu)} \ \ Changyou Wang\footnote{Department of Mathematics, Purdue University, West Lafayette, IN 47907 (wang2482@purdue.edu)} \ \ Junao Yu\footnote{University of Science and Technology of China, Hefei 230026, Anhui, China (yujunao@mail.ustc.edu.cn)}}
\title{On forward self-similar heat flow of harmonic maps}

\date{}

\begin{document}
\maketitle

\begin{abstract}
For any $k$-dimensional smooth, compact Riemannian manifold $(N, h)\subset\mathbb R^L$ without boundary,
there exists an $\varepsilon_0>0$ such that for any homogeneous of degree zero map 
$u_0(x)=\phi_0(\frac{x}{|x|}):\mathbb R^n\to N$ ($n\ge 2$), if 
$\|\nabla\phi_0\|_{L^n(\mathbb S^{n-1})}\le\varepsilon_0$ then there is a unique solution $u:\mathbb R^n\times (0,\infty)\to N$
to  the heat flow of harmonic map \eqref{HF1} and \eqref{IC}, which is forward self-similar and
belongs to $C^\infty(\R^n\times (0,\infty))\cap C^{\frac1{n}}(\R^n\times [0,\infty)\setminus \{(0,0)\})$. 

\end{abstract}

%\medskip

{\bf AMS-Classification (2020): 35Q35, 35D30, 76A15}  \\

\section{Introduction}
Let $(M, g)$ be a $n$-dimensional smooth, complete Riemannian manifold without boundary, and
$(N, h)$ be a $k$-dimensional smooth, compact Riemannian manifold  without boundary, which is isometrically
embedded into an Euclidean space $\R^L$ for $L>1$.
For $0<T\le\infty$,  a map $u\in C^2(M\times (0, T), N)$ is a heat flow of harmonic map, if it solves the heat equation of harmonic 
map:
\begin{equation}\label{HF1}
\partial_t u-\Delta_g u=\sum_{\alpha,\beta=1}^n g^{\alpha\beta}A(u)\big(\frac{\partial u}{\partial x_\alpha},\frac{\partial u}{\partial x_\beta}\big) 
\  \ {\rm{in}}\ \ M\times (0,T),
\end{equation}
where $\Delta_g$ is the Laplace-Beltrami operator on $(M,g)$ and
$A(\cdot)(\cdot, \cdot)$ denotes the second fundamental form of $N\subset\R^L$. 
Given a map $u_0\in C^2(M, N)$, we will study \eqref{HF1} together with the initial condition:
\begin{equation}
u(\cdot, 0) =u_0 \  \ {\rm{on}}\ \ M. \label{IC}
\end{equation}

The study of heat flow of harmonic maps \eqref{HF1}-\eqref{IC} was initiated by Eells and Sampson in their seminal 
paper \cite{Eells-Sampson1964}, where they proved that \eqref{HF1}-\eqref{IC} admits a unique
global smooth solution $u$ that sub-converges to a smooth harmonic map $u_\infty\in C^\infty(M, N)$
as $t\to\infty$, when the sectional curvature
of $N$ is nonpositive. For any arbitrary $(N,h)$ without curvature assumption, Chen-Struwe\cite{ChenStruwe1989}
constructed a global weak solution of \eqref{HF1}-\eqref{IC} that is smooth away from a closed set of locally finite
$n$-dimensional parabolic Hausdorff measure (cf.  Chen-Lin\cite{ChenLin1993} for $\partial M\not=\emptyset$) for $n\ge 3$,
and Struwe\cite{Struwe1985} established a unique global solution of \eqref{HF1}-\eqref{IC} with finitely many singular points
(cf. Chang\cite{Chang1989} for $\partial M\not=\emptyset$) for $n=2$. Furthermore, the existence of finite time singularity 
of \eqref{HF1}-\eqref{IC} was established by Coron-Ghidaglia\cite{CoronGhidaglia1989} and Chen-Ding \cite{ChenDing1990}
for $n\ge 3$, and by Chang-Ding-Ye\cite{ChangDingYe1992} first and Davila-del Pino-Wei\cite{DdW2020} recently for $n=2$. 
It was analyzed by Lin-Wang\cite{LinWang1999, LinWang2002} that the formation of a finite time singularity $(x_0, T_0)$  
from a locally smooth solution $u$ of \eqref{HF1} in dimensions $n\ge 3$ is closely related to the existence of backward self-similar, nontrivial solutions 
$\displaystyle\phi(\frac{x}{\sqrt{-t}})$ or steady solutions $\displaystyle\phi(\frac{x}{|x|})$, called as quasi-harmonic spheres or harmonic spheres. 
We refer interested readers to the monograph\cite{LinWang2008} and references therein for more discussions on \eqref{HF1}.

Certain equivariant backward self-similar, nontrivial solutions of \eqref{HF1} have been constructed by Fan\cite{Fan1999} for $3\le n\le 6$
and Gastel\cite{Gastel2002} for all $n\ge 3$;  while the
non-existence of backward self-similar nontrivial solutions of \eqref{HF1} has been shown under the assumption 
$\nabla u\in L^\infty([0, T], L^n)$ for $n\ge 4$ by Wang\cite{Wang2008}, see also Bizo\'n-Wasserman\cite{BW2015} for some related result. 
Typically, the construction of a continuation of the heat flow of harmonic maps \eqref{HF1}
after a singularity uses an asymptotically self-similar expander. There have been quite a few results
on forward self-similar solutions of \eqref{HF1}. Coron\cite{Coron1990} first proved that there exists a homogeneous of degree zero map
$u_0:\R^3\to\mathbb S^2$ that supports infinitely many solutions of \eqref{HF1}. Biernat-Bizo\'n\cite{BB2011} made interesting numerical
studies on both self-similar shrinker and expander of \eqref{HF1}. Germain-Rupflin\cite{GermainRupflin2011} established
the existence of forward self-similar solutions of \eqref{HF1} in the equivariant setting, and analyzed their uniqueness and stability.
Recently,  Deruelle-Lamm\cite{DL2021} established the existence of forward self-similar solutions to \eqref{HF1} in $\R^n$ for $n\ge 3$, when
the initial data $\displaystyle u_0(x)=\phi_0(\frac{x}{|x|})$ and $\phi_0\in {\rm{Lip}}(\mathbb S^{n-1}, N)$ is homotopic to 
a constant, which is smooth away from 
a compact set with finite $n$-dimensional parabolic Hausdorff measure. 
We would like to point out that, in a different but closely related context, Jia-Sverak\cite{JiaSverak2014} constructed a forward self-similar solution $u$ to the incompressible Navier-Stokes equation in $\R^3$ for any initial data $\displaystyle u_0(x)=|x|^{-1}\phi_0(\frac{x}{|x|})$, with $\phi_0\in C^\gamma(\mathbb S^2)$, such that $u\in C^\infty(\R^3\times (0,\infty))\cap C^\gamma(\R^3\times [0,\infty)\setminus\{(0,0)\})$.

\smallskip
From now on, we will assume $(M, g)=(\R^n, dx^2)$. 
A map $u:\mathbb R^n\times (0,\infty) \to N$ is called a forward self-similar heat flow of harmonic map, if 
it solves \eqref{HF1} and satisfies
$$
u(\lambda x, \lambda^2t)=u(x, t),  \ \forall (x,t)\in\R^n\times (0, \infty), \ \forall \lambda>0.
$$
This is equivalent to $\displaystyle u(x,t)=u(\frac{x}{\sqrt{t}}, 1)$ for $(x,t)\in \R^n\times (0,\infty)$.

\smallskip
Partly inspired by the main theorem by \cite{JiaSverak2014}, we establish our main result.

\begin{theorem}\label{FSSHF} For $n\ge 2$, assume $u_0\in  W^{1,n}_{\rm{loc}}(\R^n\setminus\{0\}, N)$ is homogeneous of degree zero. 
There exists an $\varepsilon_0>0$, depending on $n$ and $N$, such that 
if 
\begin{equation}\label{small0}
\|\nabla_{\mathbb S^{n-1}}u_0\|_{L^n(\mathbb S^{n-1})}\le \varepsilon_0,
\end{equation}
then there exists 
a unique solution $u\in C^\infty(\R^n\times (0, \infty), N)\cap C^{\frac1n}(\R^n\times [0, \infty)\setminus\{(0,0)\}, N)$ of the heat flow of harmonic map \eqref{HF1} and \eqref{IC}, that is forward self-similar. 
%i.e., for any $\lambda>0$,  $u(\lambda x,\lambda^2 t)=u(x,t)$ for $(x,t)\in \R^n\times (0,\infty)$. 
\end{theorem} 

First, we would like to point out the following remark.
\begin{remark}{\rm 
It follows from $u_0\in W^{1,n}_{\rm{loc}}(\R^n\setminus\{0\})$ is homogeneous of degree zero that
$$\displaystyle\|\nabla u_0\|_{L^{n,\infty}(\R^n)}\approx \displaystyle\|\nabla_{\mathbb S^{n-1}} u_0\|_{L^{n}(\mathbb S^{n-1})}\le \varepsilon_0,
$$
where $L^{n,\infty}$ denotes the weak $L^n$-space. By the Poincar\'e inequality, $u_0$ has small BMO norm:
$$
\big[u_0\big]_{\rm{BMO}(\R^n)}\lesssim \|\nabla u_0\|_{L^{n,\infty}(\R^n)}\le \varepsilon_0.
$$
Hence by \cite{Wang2011} there exists a unique solution $u\in C^\infty(\R^n\times (0,\infty), N)$ of \eqref{HF1} and \eqref{IC}. 
However, the argument by \cite{Wang2011} is insufficient for showing $u\in C^{\frac1n}(\R^n\times [0,\infty)\setminus\{(0,0)\})$,
which is the new contribution by Theorem \ref{FSSHF}.   To achieve this, we need to choose a function space $X$ that is different
from that employed by \cite{Wang2011}. It is a very interesting question to ask whether Theorem \ref{FSSHF} holds 
without the condition \eqref{small0}.}
\end{remark} 
Here we sketch the key steps of proof of Theorem \ref{FSSHF}, which are 
different from \cite{DL2021}. Instead of solving the equation of quasi-harmonic maps in $\R^n$ with
$u_0$ as the asymptotic boundary value near infinity, we will construct a mild solution $u$ of \eqref{HF1} and \eqref{IC}
for any initial data $u_0$, which has small $\|\nabla u_{0}\|_{L^{n,\infty}(\R^n)}$.
More precisely, 
\begin{enumerate}
\item[Step 1.] Seek a mild solution of \eqref{HF1}-\eqref{IC} $u(\cdot, t)=\widehat{u_0}(t)+v(\cdot, t)$, 
where $\widehat{u_0}(t)=e^{t\Delta} u_0$ and $v(x,t)$ is a fixed point of the contractive solution operator $\mathcal{S}:\mathbb B_{\delta}(0)
\to \mathbb B_{\delta}(0)$ given by 
$$
\mathcal{S}(v)(\cdot, t):=\int_0^t e^{(t-s)\Delta}A(\widehat{u_0}+v)(\nabla\widehat{u_0}+\nabla v, \nabla\widehat{u_0}+\nabla v)(s)\,ds, \ t>0, 
$$  
where $\mathbb B_\delta(0)\subset X$ is the $\delta$-ball center at $0$, 
and  $X$ is the Banach space with norm
$\|v\|_{X}:=\|v\|_{L^\infty(\R^n\times \R_+)}+\sup_{t\ge 0}t^\frac14\|\nabla v\|_{L^{2n}(\R^n)}$.
Here the assumption on $u_0$ guarantees  
\begin{align*}
\sup_{t\ge 0}
t^{\frac14}\big\|\nabla\widehat{u_{0}}\big\|_{L^{2n}(\R^n)} 
\le C\|\nabla u_{0}\|_{L^{n,\infty}(R^n)}\le C\delta.
\end{align*}
Moreover, 
\begin{align}\label{smallness1}
\sup_{t\ge 0} \|\nabla u(t)\|_{L^{n,\infty}(\R^n)}\le C\delta,
\end{align}
and
\begin{align}\label{smallness2}
\sup_{x\in\R^n, R>0} R^{-n}\int_{B_R(x)\times [0, R^2]} |\nabla u|^2\le C\delta^2.
\end{align}
%Furthermore,  if $u_0$ is homogeneous of degree zero, then $u$ can be shown to be forward self-similar, 
%i.e. $u(x,t)=U\big(\frac{x}{\sqrt{t}})$, and $U\in C^\infty(\R^n)$.
\item[Step 2.]  \eqref{smallness1} ensures the smallness of BMO norm of $u(t)$ 
so that 
\begin{equation}\label{rigid ineq}
\big\|\Delta u\big\|_{L^2(B_R(x))}\le C\Big(R^{-1}\|\nabla u\|_{L^2(B_R(x))}+\|\Delta u+A(u)(\nabla u, \nabla u)\|_{L^2(B_R(x))}\Big)
\end{equation}
holds for all $B_R(x)\subset\R^n$. 
This and the smallness of renormalized energy \eqref{smallness2} enable us to employ a blowing up argument to show 
the renormalized energy enjoys
a decay property that yields the H\"older continuity of $u$ at any $(x_0, 0)$, with $x_0\not=0$, if in addition
$u_0$ is homogeneous of degree zero. 
\end{enumerate}

\section{Proof of Theorem \ref{FSSHF}} 

This section is devoted to the proof of Theorem \ref{FSSHF}.  First, we will apply the fixed point argument to obtain a mild solution $u$ 
of \eqref{HF1} for any initial data $u_0$,
with small $\|\nabla u_{0}\|_{L^{n,\infty}(\R^n)}$, 
that has sufficiently high Sobolev integrability and small $\|\nabla u\|_{L^{n,\infty}(\R^n)}$. 
%Then we will show such a constructed solution is smooth in $\R^n\times (0,\infty)$
%and is H\"older continuous up to $t=0$ with exponent $\gamma$,  except at $(0,0)$. 

%First we recall the following approximation property of the initial data $u_0:\R^n\to N$.

%\begin{lemma}\label{3infty-approx} {For any $\delta>0$, there exist $M_0>0$ depending on $\delta$ and $\|\nabla u_0\|_{L^{n,\infty}(\R^n)}$, 
%$M_1>0$ depending on $N$, and $u_{01}, u_{02}:\R^n\to \R^L$  such that 
%\begin{equation}
%\max\big\{\|u_{01}\|_{L^\infty(\R^n)},\ \|u_{02}\|_{L^\infty(\R^n)}\big\}\le M_1,
%\end{equation}
%and
%\begin{equation}\label{approx_id}
%u_0=u_{01}+u_{02}\ \ {\rm{in}}\ \ \R^n,
%\end{equation}
%\begin{equation}\label{3infty-est}
%\|\nabla u_{01}\|_{L^{n, \infty}(\R^n)}\le\delta,  \ \  \|\nabla u_{02}\|_{L^{2n}(\R^n)}\le M_0.
%\end{equation}}
%\end{lemma}

\bigskip
For $t\ge 0$, let $\widehat{u_{0}}(t)=e^{t\Delta} u_{0}: \R^n\to\R^L$ denote the Caloric extension of $u_{0}$. Namely,
$$(\partial_t-\Delta)\widehat{u_{0}}=0 \ \ {\rm{in}}\ \ \R^n\times (0,\infty),  \ \ \widehat{u_{0}}(0)=u_{0} \ \ {\rm{on}}\ \ \R^n.$$
By the standard Lorentz space $L^{p, q}$-estimate on convolution operators through the heat kernel, we obtain
\begin{lemma}\label{conv-est} If $u_0:\R^n\to N$ satisfies $\nabla u_0\in L^{n,\infty}(\R^n)$, then, for any $t> 0$,
\begin{equation}\label{3infty-est1}
\|\nabla \widehat{u_{0}}(t)\|_{L^{n,\infty}(\R^n)}\le \|\nabla {u_{0}}\|_{L^{n,\infty}(\R^n)},
\end{equation}
\begin{equation}\label{2n-est2}
t^{\frac14}\|\nabla \widehat{u_{0}}(t)\|_{L^{2n}(\R^n)}\le C\|\nabla {u_{0}}\|_{L^{n,\infty}(\R^n)},
\end{equation}
\begin{equation}\label{p-est3}
t^{\frac{n}2(\frac{1}{n}-\frac{1}{p})}\|\nabla \widehat{u_{0}}(t)\|_{L^{p}(\R^n)}\le C\|\nabla {u_{0}}\|_{L^{n,\infty}(\R^n)}, 
\ \forall n<p\le\infty.
\end{equation}
\end{lemma}

\bigskip
\pf We leave the proof to interested readers.
\qed

\bigskip
Now we are ready to establish the following  theorem on \eqref{HF1} and \eqref{IC}.

\begin{theorem}\label{FSSHF1} There exist a small $\delta_0>0$ such that if $\|\nabla u_0\|_{L^{n,\infty}(\R^n)}\le\delta_0$,
then there is a unique solution $u:\mathbb R^n\times [0, \infty)\to N$ of \eqref{HF1} and \eqref{IC}  such that 
$v:=u-\widehat{u_{0}}$ and $u$ satisfy the following estimates:
\begin{align}\label{xt-est-v1}
&\big\|v\big\|_{L^\infty(\R^n\times \R_+)}+\sup_{t\ge 0} t^{\frac14}\big\|\nabla v(t)\big\|_{L^{2n}(\R^n)} \le C\delta_0, \\
&\big\|\nabla v(t)\big\|_{L^{n}(\R^n)} \le C\delta^2_0, \ \forall t>0,  \label{xt-est-v2}
\end{align}
and
\begin{align}
&\big\|\nabla u(t)\big\|_{L^{n,\infty}(\R^n)}\le C\delta_0, \ \forall t>0,
\label{ninfty-est-u} \\
&R^{-2}\int_{B_R(x)\times [0, R^2]}|\nabla u|^n \le C\delta_0^2, \ \forall x\in\R^n, \ \forall R>0. \label{normalized-n-est-u}
\end{align}
If, in addition, $u_0$ is homogeneous of degree zero, then $u$ is forward self-similar.
\end{theorem}

\bigskip
\pf
Define the Banach space $X$ by 
\begin{align*}
X&:=
\Big\{v:\mathbb R^n\times [0, \infty)\to \R^L: \  v(\cdot, 0)=0 \ {\rm{on}}\ \R^n, \\
&\qquad \|v\|_{X}:=\|v\|_{L^\infty(\R^n\times [0,\infty))}+\sup_{t\ge 0} t^{\frac14}\big\|\nabla v\big\|_{L^{2n}(\R^n)}<\infty\Big\}.
\end{align*}
For $\delta\ge \delta_0$, define the $\delta$-ball, with center $0$, in $X$ by 
\begin{equation}\label{delta-ball}
\mathbb{B}_\delta(0)=\Big\{ v\in X: \ \|v\|_{X}\le \delta\Big\}.
\end{equation}

Recall that, since $N\subset\mathbb R^L$ is compact and smooth, there exists a $\delta_N>0$ such that the nearest point
projection map  $\Pi_N$ from the $\delta_N$-neighborhood, $N_{\delta_N}=\big\{y\in\R^L: {\rm{dist}}(y, N)<\delta_N\big\}$, to $N$
is smooth. The second fundamental form of $N$ is defined by 
$$A(y)(X, Y)=D_y^2\Pi_N(y)(X, Y), \ \forall y\in N, \ X, Y\in T_yN.$$
Now, let $\widehat{A}(\cdot)$ be a smooth extension of the second fundamental form $A$ of $N$ such that
\begin{align*}
\widehat{A}(y)=A(\Pi_N(y))  \ \forall y\in N_{\delta_N}, \
\displaystyle\sup_{y\in\R^L}|\widehat{A}(y)|\le C(N),
\end{align*}
where $C(N)>0$ is a constant depending only on $N$.

Next, define the solution operator $\mathcal{S}: \mathbb{B}_\delta(0)\to X$ by letting $w=\mathcal{S} v$ and 
\begin{equation}
w(\cdot, t)=\int_0^t e^{(t-s)\Delta} \widehat{A}(u)(\nabla u,\nabla u)(s)\,ds, \ t\ge 0,
\end{equation}
where $u(s)=\widehat{u_{0}}(s)+v(s)$ for $s\ge 0$.  It is easy to see that
$w(x,0)=0$ in $\R^n$.  Applying H\"older's inequality, we can estimate that for any $t>0$, 
\begin{align}\label{l-infinity-est1}
\|w(t)\|_{L^\infty(\R^n)}&\le C \int_0^t \|{(t-s)^{-\frac{n}2}}e^{-\frac{|y|^2}{4(t-s)}}\|_{L^{\frac{n}{n-1}}(\R^n)}
\||\nabla u(s)|^2\|_{L^{n}(\R^n)}\,ds\nonumber\\
&\le C\big(\int_{0}^t (t-s)^{-\frac12} s^{-\frac12}\,ds\big)\big(\|\widehat{u_{0}}\|_{X}^2
+\|v\|_{X}^2\big)\le C(\delta_0^2+\delta^2)\le \frac{\delta}{2},
\end{align}
provided $\delta_0\le \delta$ and $\delta>0$ is chosen sufficiently small. 

Also, by applying the inequality
\begin{equation}\label{lp-est}
\|\nabla e^{t\Delta} f\|_{L^{2n}(\R^n)} \le Ct^{-\frac14-\frac{n}{2p}} \|f\|_{L^p(\R^n)}, \ \forall 1\le p\le 2n,
\end{equation}
we can estimate that for any $t>0$,
\begin{align*}
\|\nabla w(t)\|_{L^{2n}(\R^n)} 
&\le C \int_0^t (t-s)^{-\frac34}\|\nabla u(s)\|_{L^{2n}(\R^n)}^2\,ds\\
&\le C \int_0^t (t-s)^{-\frac34}\big(\|\nabla\widehat{u_{0}}\|_{L^{2n}(\R^n)}^2
+\|\nabla v\|_{L^{2n}(\R^n)}^2\big)(s)\,ds\\
&\le C \big(\int_0^t (t-s)^{-\frac34}s^{-\frac12}\,ds\big) \|\nabla{u_{0}}\|_{L^{n,\infty}(\R^n)}^2\\
&\ \ +C\big(\int_0^t (t-s)^{-\frac34} s^{-\frac12} \,ds\big) \|v\|_{X}^2\\
&\le C t^{-\frac14} \big(\|v\|_{X}^2+\delta_0^2).
\end{align*}
Thus 
\begin{equation*}
\sup_{t\ge 0}t^\frac14\|\nabla w\|_{L^{2n}(\R^n)}\le C\big(\|v\|_{X}^2+\delta_0^2)
<\frac{\delta}2,
\end{equation*} 
provided $\delta>0$ is chosen sufficiently small. Hence 
$\mathcal{S}(\mathbb{B}_{\delta}(0))\subset \mathbb B_{\delta}(0)$. 

\medskip
Next we want to show that $\mathcal{S}$ is a contraction map on $\mathbb B_{\delta}(0)$. 
For this, let $v_1, v_2\in \mathbb B_{\delta}(0)$, and
$w_1=\mathcal{S}(v_1)$ and $w_2=\mathcal{S}(v_2)$. Denote $u_i=\widehat{u_{0}} +v_i$ for $i=1,2$. We obtain
\begin{align*}
&|(w_1-w_2)(t)|=\big|\int_0^t e^{(t-s)\Delta} \big(\widehat{A}(u_1)(\nabla u_1,\nabla u_1)
-\widehat{A}(u_2)(\nabla u_2,\nabla u_2)\big)(s)\,ds\big|\\
&\le C \int_0^t e^{(t-s)\Delta}\Big(|v_1-v_2||\nabla u_1|^2+|\nabla (v_1-v_2)|(|\nabla \widehat{u_{0}}|
+|\nabla v_1|+|\nabla v_2|)\Big)(s)\,ds.
\end{align*}
Hence, similar to \eqref{l-infinity-est1}, we can estimate that for $t\ge 0$, 
\begin{align*}
&\|(w_1-w_2)(t)\|_{L^\infty(\R^n)}\\
&\le C\|v_1-v_2\|_{X} \big(\int_0^t (t-s)^{-\frac12} s^{-\frac12}\,ds\big) (\|\widehat{u_{0}}\|_{X}^2+\|v_1\|_{X}^2)\\
&+C\|v_1-v_2\|_{X}\big(\int_0^t (t-s)^{-\frac12} s^{-\frac12}\,ds\big)
\big(\|\widehat{u_{0}}\|_{X_T}+\|v_1\|_{X}+\|v_2\|_{X}\big)\\
&\le C\big(\|\widehat{u_{0}}\|_{X}^2+\|v_1\|_{X}^2+\|\widehat{u_{0}}\|_{X}+\|v_1\|_{X}+\|v_2\|_{X}\big)\|v_1-v_2\|_{X}\\
&\le C\delta \|v_1-v_2\|_{X}
\end{align*}
so that
\begin{align}\label{xt-est1}
 \|w_1-w_2\|_{L^{\infty}(\R^n\times \R_+)}
\le C\delta \|v_1-v_2\|_{X}.
\end{align}
Furthermore, 
\begin{align*}
&\|\nabla(w_1-w_2)(t)\|_{L^{2n}(\R^n)}\\
&\le C\|v_1-v_2\|_{X} \big(\int_0^t(t-s)^{-\frac34} s^{-\frac12}\,ds\big)(\|\widehat{u_{0}}\|_{X}^2+\|v_1\|_{X}^2)\\
&+C\|v_1-v_2\|_{X}\big(\int_0^t (t-s)^{-\frac34} s^{-\frac12}\,ds\big)
\big(\|\widehat{u_{0}}\|_{X}+\|v_1\|_{X}+\|v_2\|_{X}\big)\\
&\le C\|v_1-v_2\|_{X}t^{-\frac14}\big(\|\widehat{u_{0}}\|_{X}^2+\|v_1\|_{X}^2\big)\\
&+C\|v_1-v_2\|_{X}t^{-\frac14}
\big(\|\widehat{u_{0}}\|_{X}+\|v_1\|_{X}+\|v_2\|_{X}\big)
\end{align*} 
so that
\begin{align}\label{xt-est2}
\sup_{t\ge 0} t^\frac14 \|\nabla(w_1-w_2)(t)\|_{L^{2n}(\R^n)}
\le C\delta \|v_1-v_2\|_{X}.
\end{align}
Combining \eqref{xt-est1} with \eqref{xt-est2}, we obtain 
\begin{align}\label{xt-est3}
\|w_1-w_2\|_{X}\le \theta \|v_1-v_2\|_{X},
\end{align}
where $0<\theta=C\delta<1$, provided $\delta>0$ is chosen sufficiently small.
Therefore we can apply the fixed point theorem to deduce that there exists a unique
$v\in \mathbb{B}_\delta(0) \subset X$ such that 
\begin{align}\label{sol1}
v(t)=\int_0^t e^{(t-s)\Delta}\widehat{A}(\widehat{u_0}+v)(\nabla\widehat{u_0}+\nabla v, \nabla\widehat{u_0}+\nabla v)(s)\,ds,
\ t>0.
\end{align}
We can apply 
\begin{equation}\label{lp-est1} 
\|\nabla e^{t\Delta} f\|_{L^n(\R^n)}\le Ct^{-\frac{n}{2p}} \|f\|_{L^p(\R^n)}, \ 1\le p\le n,
\end{equation}
to estimate $\|\nabla v\|_{L^n(\R^n)}$ by
\begin{align}\label{Ln-est-v}
\|\nabla v(t)\|_{L^n(\R^n)}
&\le C\int_0^t (t-s)^{-\frac12} \|\nabla u(s)\|_{L^{2n}(\R^n)}^2\,ds\nonumber\\
&\le C\Big(\|\widehat{u_{0}}\|_{X}^2+\|v\|_{X}^2\Big) \int_0^t (t-s)^{-\frac12} s^{-\frac12}\,ds \le C\delta^2.
\end{align}

If we define $u(x,t)=\widehat{u_0}(x,t)+v(x,t)$ for $(x,t)\in \R^n\times \R_+$, then $\nabla u(t)\in L^{n,\infty}(\R^n)$, and
\begin{align}\label{ninfty-est}
\|\nabla u(t)\|_{L^{n,\infty}(\R^n)}\le \|\nabla \widehat{u_{0}}(t)\|_{L^{n,\infty}(\R^n)}
+\|\nabla v(t)\|_{L^{n}(\R^n)}\le C\delta.
\end{align}

For $K>0$ and $(x,t)\in\R^n\times (0, \infty)$, define
$$c_{0}=\frac{1}{|B_K|}\int_{B_K} u_{0}(x-\sqrt{t}y)\,dy.$$
Then 
\begin{align*}
&|\widehat{u_0}(x,t)-c_0|
\le\big |\int_{\R^n} \frac{1}{(4\pi)^{\frac{n}2}}e^{-\frac{|y|^2}4}(u_{0}(x-\sqrt{t}y)-c_{0})\,dy\big|\\
&\le \int_{B_K}|u_{0}(x-\sqrt{t}y)-c_{0}|\,dy+\int_{\R^n\setminus B_K} e^{-\frac{|y|^2}4}|u_{0}(x-\sqrt{t}y)-c_{0}|\,dy\\
&\le K^n \|\nabla u_{0} \|_{L^{n,\infty}(B_{K\sqrt{t}})}+C_N\int_K^\infty e^{-\frac{r^2}4} r^{n-1}\,dr\\
&\le K^n \|\nabla u_{0}\|_{L^{n,\infty}(\R^n)}+o_K(1)\le K^n\delta+o_K(1),
\end{align*}
where $\displaystyle\lim_{K\to \infty} o_K(1)=0.$

On the other hand, since $u_0(\R^n)\subset N$, it follows from the Poincar\'e inequality that
\begin{align*}
{\rm{dist}}(c_0, N)&\le \frac{1}{|B_K|}\int_{B_K}|u_0(x-\sqrt{t}y)-c_0|\,dy
\le C\|\nabla u_{0}\|_{L^{n,\infty}(B_{K\sqrt{t}})}\le C\delta.
\end{align*}
Therefore we obtain that
\begin{align}\label{distance}
&{\rm{dist}}(u(x,t), N)\le {\rm{dist}}(\widehat{u_0}(x,t), N)+|v(x,t)|\nonumber\\
&\le |\widehat{u_0}(x,t)-c_0|+{\rm{dist}}(c_0, N)+|v(x,t)|\le(1+ K^n)\delta+o_K(1)<\delta_N,
\end{align} 
provided $K$ is chosen sufficiently large and $\delta$ is chosen sufficiently small.

It follows from \eqref{sol1}  that $u$ solves
\begin{equation}\label{HF2}
\begin{cases}
\partial_t u-\Delta u=\widehat{A}(u)(\nabla u, \nabla u) &\ \ {\rm{in}}\ \ \R^n\times (0, \infty),\\
u\big|_{t=0}=u_0&\ \ {\rm{on}}\ \ \R^n.
\end{cases}
\end{equation}
With the help of \eqref{distance},  one can argue as in the proof of Theorem 1.3 of \cite{Wang2011} to 
show that $u(\R^n\times \R_+)\subset N$ and hence
$\widehat{A}(u)(\nabla u, \nabla u)=A(u)(\nabla u, \nabla u)$. In particular, $u$ solves the equation \eqref{HF1}
for the heat flow of harmonic map.

Applying \eqref{lp-est} to $\nabla\widehat{u_{0}}$, we obtain that, for $t>0$, 
\begin{equation}\label{lp-est2}
\|\nabla\widehat{u_{0}}(t)\|_{L^p(\R^n)} \le Ct^{-\frac12+\frac{n}{2p}}\|\nabla u_{0}\|_{L^{n,\infty}(\R^n)}
\le C\delta t^{-\frac12+\frac{n}{2p}},\ \forall p>n.
\end{equation}
Thus we obtain that for any $x_0\in \R^n$ and $R>0$, 
\begin{align}
&R^{-2}\int^{R^2}_{0}\int_{B_R(x_0)} |\nabla u|^n
\le CR^{-2}\int^{R^2}_{0}\int_{B_R(x_0)} \big(|\nabla \widehat{u_{0}}|^n+|\nabla v|^n\big)\nonumber\\
&\le C\sup_{0<t\le R^2} \|\nabla v(t)\|_{L^n(\R^n)}^n\nonumber\\
&\quad+CR^{-2+n(1-\frac{n}{p})}\int^{R^2}
_{0}\big(\int_{B_R(x_0)} |\nabla \widehat{u_{0}}(t)|^p\big)^{\frac{n}{p}} \,dt\ \ \ \ \Big(n<p<\begin{cases} \frac{n^2}{n-2} & n\ge 3\\ \infty & n=2\end{cases}\Big) \nonumber\\
&\le C\sup_{0<t\le R^2} \|\nabla v(t)\|_{L^n(\R^n)}^n+CR^{-2+n(1-\frac{n}{p})}\int^{R^2}_{0}t^{(-\frac12+\frac{n}{2p})n} \,dt \|\nabla u_{0}\|_{L^{n,\infty}(\R^n)}^n\nonumber\\
&\le C\delta^n.
\end{align}

Since $u_0(x)=u_0(\frac{x}{|x|}), x\not=0$, is homogeneous of degree zero, it follows that 
$\widehat{u_0}(t)=e^{t\Delta} u_0$ is self-similar, i.e.,
$$
\widehat{u_0}(\lambda x, \lambda^2 t)=\widehat{u_0}(x, t), \ \forall (x,t)\in \R^n\times (0,\infty), \ \forall \lambda>0. 
$$
Thus for any $\lambda>0$, it follows from \eqref{sol1} that
$v_\lambda(x,t)=v(\lambda x, \lambda^2 t), (x,t)\in \R^n\times [0,\infty)$, also satisfies
\begin{align}\label{sol2}
v_\lambda(t)=\int_0^t e^{(t-s)\Delta}\widehat{A}(\widehat{u_0}+v_\lambda)(\nabla\widehat{u_0}+\nabla v_\lambda, \nabla\widehat{u_0}+\nabla v_\lambda)(s)\,ds,
\ t>0.
\end{align}
It is easy to verify $\|v_\lambda\|_{X}\le \|v\|_{X}$
so that $v_\lambda\in \mathbb{B}_\delta(0)\subset X$ and hence
$v_\lambda\equiv v, \ \forall \lambda>0.$
\qed

\bigskip
Next we will show that $u$ by Theorem \ref{FSSHF1} enjoys regularity in Theorem \ref{FSSHF}.

\begin{proposition}\label{interior-smooth} Assume $u_0\in W^{1,n}_{\rm{loc}}(\R^n\setminus\{0\}, N)$
is homogeneous of degree zero, and $u:\mathbb R^n\times \R_+\to N$ is the solution of \eqref{HF1} by Theorem \ref{FSSHF1}.
Then $u\in C^\infty(\R^n\times (0, \infty), N)$. 
\end{proposition} 
\pf By Theorem \ref{FSSHF1}, $u$ is a forward self-similar solution of \eqref{HF1} on $\R^n\times (0,\infty)$ and can be written
as $\displaystyle u(x,t)={U}(\frac{x}{\sqrt{t}})$, where $U:\R^n\to N$ solves the following equation for quasi-harmonic maps: 
\begin{equation}\label{qhm}
\Delta U+\frac12y\cdot\nabla U+A(U)(\nabla U,\nabla U)=0 \ \ {\rm{in}}\ \ \R^n.
\end{equation} 
Moreover, it follows from Theorem \ref{FSSHF1} that
\begin{equation}\label{2n-bound}
\|\nabla U\|_{L^{2n}(\R^n)}\equiv\sup_{t>0} t^{\frac14}\|\nabla u(t)\|_{L^{2n}(\R^n)} \le \delta.
\end{equation}
From this, we can apply the $W^{2,n}$-estimate to \eqref{qhm} to conclude that for any $R>0$, $U\in W^{2,n}(B_R)$,
and 
\begin{equation}\label{2n-bound1}
\|\nabla^2 U\|_{L^n(B_R)} \le C(n, R)(R^3+\|\nabla U\|_{L^{2n}(B_{2R})}^2) \le C(n, R)(\delta^2+R^3). 
\end{equation}
This, combined with the Sobolev embedding inequality, can further imply that $\nabla U\in L^{q}(B_R)$ for any $1<q<\infty$, and
\begin{align}\label{2n-bound2}
\|\nabla U\|_{L^q(B_R)} \le C(q,R)\big(\|\nabla U\|_{L^n(B_R)}+R^3+\|\nabla U\|_{L^{2n}(B_{2R})}^2\big)
\le C(q, R)(\delta^2+R^3). 
\end{align}
This, combined with $W^{2,q}$-estimate,  further implies that $U\in W^{2,q}(B_R)$ for any $q\in (1,\infty)$, and
\begin{align}\label{2n-bound3}
\|\nabla^2 U\|_{L^q(B_R)} \le C(q,R)(\delta^2+R^3). 
\end{align}
By Morrey's embedding theorem, we conclude that $U\in C^{1,\alpha}(B_R)$ for any $\alpha\in (0,1)$. 
Applying the Schauder theory, one can show
that $U\in C^\infty(\R^n)$ and 
\begin{equation}\label{ck-bound}
\big[U\big]_{C^k(B_R)}\le C(k, R) (\delta^2+R^3), \ \ \forall k\ge 1 \ {\rm{and}}\ R\ge 1.
\end{equation} 
This completes the proof. \qed

\begin{theorem}\label{initial-Holder} Assume $u_0\in W^{1,n}_{\rm{loc}}(\R^n\setminus\{0\}, N)$
is homogeneous of degree zero, and $u:\mathbb R^n\times [0,\infty)\to N$ is the solution of \eqref{HF1} constructed by Theorem \ref{FSSHF1}.
Then $u\in C^{\frac1n}(\R^n\times [0, T]\setminus\{(0,0)\}, N)$.
\end{theorem} 

\bigskip
The proof of Theorem \ref{initial-Holder} is divided into several Lemmas. 

\begin{lemma} \label{LEI1} For any nonnegative $\eta\in C^\infty_0(\R^n)$, it holds that
\begin{align}\label{LEI10} 
&\int_{\R^n}|\nabla u(t)|^2\eta^2+\int_{\R^n\times [0,t]} |\Delta u+A(u)(\nabla u,\nabla u)|^2\eta^2\nonumber\\
&\le \int_{\R^n}|\nabla  u_0|^2\eta^2+4\int_{\R^n\times [0,t]} |\nabla u|^2|\nabla\eta|^2.
\end{align}
\end{lemma}

\pf From Proposition \ref{interior-smooth}, $u\in C^\infty(\R^n\times (0, \infty))$ and we can multiply \eqref{HF1}
by $(\Delta u+A(u)(\nabla u,\nabla u))\eta^2$ and integrate the resulting equation to obtain
\begin{align*}
&\frac{d}{dt}\int_{\R^n}|\nabla u|^2\eta^2 +2\int_{\R^n}|\Delta u+A(u)(\nabla u,\nabla u)|^2\eta^2\\
&=-4\int_{\R^n}\nabla u\cdot\partial_t u \eta\nabla\eta\le \int_{\R^n}|\Delta u+A(u)(\nabla u,\nabla u)|^2\eta^2
+4\int_{\R^n}|\nabla u|^2|\nabla \eta|^2.
\end{align*}
Hence we arrive at
\begin{align*}
\frac{d}{dt}\int_{\R^n}|\nabla u|^2\eta^2 +\int_{\R^n}|\Delta u+A(u)(\nabla u,\nabla u)|^2\eta^2
\le 4\int_{\R^n}|\nabla u|^2|\nabla \eta|^2.
\end{align*}
Integrating this inequality from $0$ to $t$ yields \eqref{LEI10}. \qed

\begin{corollary}\label{LEI2} Let $u:\R^n\times \R_+\to N$ be the solution given by Theorem 
\ref{FSSHF1}. Then, for any $x_0\in\R^n$ and $R>0$, it holds  that
\begin{align}\label{LEI20}
&\sup_{0\le t\le R^2}\int_{B_{\frac{R}2}(x_0)}|\nabla u(t)|^2+\int_{B_{\frac{R}2}(x_0)\times [0, R^2]} |\Delta u+A(u)(\nabla u, \nabla u)|^2\nonumber\\
&\le \int_{B_{R}(x_0)}|\nabla u_0|^2+\frac{64}{R^2}\int_{B_{R}(x_0)\times [0, R^2]} |\nabla u|^2.
\end{align}
\end{corollary}
\pf Let $\eta\in C^\infty_0(\R^n)$ be a cut-off function of $B_{\frac{R}2}(x_0)$ such that 
$$\eta\equiv 1 \ {\rm{in}}\ B_{\frac{R}2}(x_0), \ \eta\equiv 0 \ {\rm{in}}\ \R^n\setminus B_R(x_0), \ |\nabla\eta|\le \frac{4}{R}.$$
Substituting $\eta$ into \eqref{LEI10} and varying $t$ from $0$ to $R^2$ yields \eqref{LEI20}.
\qed

\bigskip
Now we are ready to establish an energy decay property of the solution $u$ under a smallness condition. 
Recall from \eqref{ninfty-est-u} and \eqref{normalized-n-est-u}
of Theorem \ref{FSSHF1} that $u$ satisfies 
\begin{equation}\label{bmo-est}
\sup_{t\ge 0} \|\nabla u(t)\|_{L^{n,\infty}(\R^n)} \le C\delta, 
\end{equation}
\begin{equation}\label{normalized-2-est-u}
R^{-n}\int_{B_R(x)\times [0, R^2]} |\nabla u|^2 \le C\delta^2, \ \forall x\in\R^n,  \ R>0.
\end{equation}

\begin{lemma}\label{energy-decay1} There exist $\varepsilon_0>0$, $0<R_0\le \frac14$, $\theta_0\in (0,1)$, and $C_0>0$ such that 
for any $x_0\in \R^{n}$, with $\frac12\le |x_0|\le 2$,  if 
\begin{equation}\label{small1}
R_0^{-n}\int_{B_{R_0}(x_0)\times [0,R_0^2]}|\nabla u|^2\le \varepsilon_0^2,
\end{equation}
then 
\begin{equation}\label{energy-decay10}
(\theta_0R_0)^{-n}\int_{B_{\theta_0R_0}(x_0)\times [0,(\theta_0R_0)^2]}|\nabla u|^2\le \frac12\max\Big\{
R_0^{-n}\int_{B_{R_0}(x_0)\times [0,R_0^2]}|\nabla u|^2, C_0 R_0^{\frac2n}\Big\}.
\end{equation}
\end{lemma} 
\pf We argue by contradiction. Suppose that the conclusion were false. Then for any fixed $\theta\in (0,1)$,
there exist $x_i\in \mathbb R^{n}$ with $\frac12\le |x_i|\le 2$, $\varepsilon_i\to 0$, $R_i\to 0$ such that
\begin{equation}\label{small2}
R_i^{-n}\int_{B_{R_i}(x_i)\times [0,R_i^2]}|\nabla u|^2=\varepsilon_i^2,
\end{equation}
while 
\begin{equation}\label{energy-decay11}
(\theta R_i)^{-n}\int_{B_{\theta R_i}(x_i)\times [0,(\theta R_i)^2]}|\nabla u|^2>\frac12\max\Big\{
\varepsilon_i^2, i R_i^{\frac2n}\Big\}.
\end{equation}
It follows from \eqref{small2} and \eqref{energy-decay11} that 
\begin{equation}\label{r-est}
R_i\le \frac{2^{\frac{n}2}\varepsilon_i^n}{i^{\frac{n}2}\theta^{\frac{n^2}2}}.
\end{equation}
Applying \eqref{bmo-est}, we have
\begin{equation}\label{bmo-est1}
\sup_{t\ge 0} \|\nabla u(t)\|_{L^{n,\infty}(\R^n)} \le C\delta.
\end{equation}

Now we define 
$$u^i(x,t)=u(x_i+R_i x, R_i^2 t), \ (x,t)\in B_1(0)\times [0,1],  \ \ u^i_0(x)=u_0(x_i+R_i x), 
\ x\in B_1(0).$$
Then $u_i$ satisfies
\begin{equation}\label{HF3}
\begin{cases}
\partial_t u^i=\Delta u^i+A(u^i)(\nabla u^i,\nabla u^i)\ \  {\rm{in}}\ \ B_1(0)\times [0, 1],\\
u^i(\cdot,0)=u^i_0 \  \ {\rm{in}}\ \  B_1(0),
\end{cases}
\end{equation}
\begin{equation}\label{small3}
\displaystyle\int_{B_1(0)\times [0,1]}|\nabla u^i|^2=\varepsilon_i^2,\ \ 
\displaystyle\theta^{-n}\int_{B_{\theta }(0)\times [0, \theta ^2]}|\nabla u^i|^2>\frac12\max\Big\{
\varepsilon_i^2, i R_i^{\frac2n}\Big\},
\end{equation}
and
\begin{equation} \label{bmo-est2}
\sup_{0\le t\le 1} \|\nabla u^i(t)\|_{L^{n,\infty}(B_{1}(0))} \le C\delta.
\end{equation}
Since $u_0$ is of homogeneous degree zero, by direct calculations we  have
\begin{align} \label{u0-est}
&\int_{B_1(0)}|\nabla u^i_0|^2=R_i^{2-n}\int_{B_{R_i}(x_i)}|\nabla u_0|^2\nonumber\\
&\le \big(\int_{B_{R_i}(x_i)}|\nabla u_0|^n\big)^{\frac2n}\le \Big(\int_{|x_i|-R_i}^{|x_i|+R_i}\frac{dr}{r} \int_{\mathbb S^{n-1}}|\nabla_{\rm{tan}} u_0|^n\,d\sigma\Big)^{\frac2n}\nonumber\\
&\le C R_i^{\frac2n} \|u_0\|_{W^{1,n}(\mathbb S^{n-1})}^{2}.
\end{align}
Hence, by applying \eqref{LEI20} to $u^i$, one has
\begin{align}\label{LEI30}
&\sup_{0\le t\le 1}\int_{B_{\frac{1}2}(0)}|\nabla u^i(t)|^2+\int_{B_{\frac{1}2}(0)\times [0, 1]} |\Delta u^i+A(u^i)(\nabla u^i, \nabla u^i)|^2\nonumber\\
&\le \int_{B_{1}(0)}|\nabla u_0^i|^2+C\int_{B_{1}(0)\times [0, 1]} |\nabla u^i|^2
\le C\big(R_i^{\frac2n}+\varepsilon_i^2\big).
\end{align}

Applying the interpolation inequality \cite{Adams-Frazier} (see also \cite{Wang2004} Proposition 4.3),  we have 
\begin{equation}\label{l4-est1}
\big\|\nabla u^i(t)\big\|_{L^4(B_\frac12(0))}^4
\le  C\big\|\nabla u^i(t)\big\|_{L^{n,\infty}(B_{\frac12}(0))}^2 \Big\||\nabla u^i(t)|+|\nabla^2 u^i(t)|\Big\|_{L^2(B_\frac12(0))}^2.
\end{equation} 
This, combined with \eqref{bmo-est2}, yields
\begin{align}\label{l4-est2}
&\int_{B_\frac12(0)\times [0,1]} |\nabla u^i|^4
\le  C\sup_{0\le t\le 1}\|\nabla u^i(t)\|_{L^{n,\infty}(B_{\frac12}(0))}^2 \int_{B_\frac12(0)\times [0,1]}
(|\nabla u^i|^2+|\nabla^2 u^i|^2)\nonumber\\
&\le C\delta\int_{B_\frac12(0)\times [0,1]}
\big(|\nabla u^i|^2+|\Delta u^i|^2\big)\nonumber\\
&\le C\delta\int_{B_\frac12(0)\times [0,1]}
\big(|\nabla u^i|^2+|\Delta u^i+A(u^i)(\nabla u^i,\nabla u^i)|^2+|\nabla u^i|^4\big).
\end{align}
Therefore, if we choose $\delta>0$ sufficiently small  so that $C\delta\le \frac12,$
then
\begin{align}\label{l4-est3}
\int_{B_\frac12(0)\times [0,1]} |\nabla u^i|^4
&\le \int_{B_\frac12(0)\times [0,1]}
\big(|\nabla u^i|^2+|\Delta u^i+A(u^i)(\nabla u^i,\nabla u^i)|^2\nonumber\\
&\le \int_{B_1(0)}|\nabla u_0^i|^2+C\int_{B_1(0)\times [0, 1]}|\nabla u^i|^2.
\end{align}
Substituting \eqref{l4-est3} into \eqref{LEI30}, one can improve \eqref{LEI30} into the following estimate:
\begin{align}\label{LEI40}
&\sup_{0\le t\le 1}\int_{B_{\frac{1}2}(0)}|\nabla u^i(t)|^2+\int_{B_{\frac{1}2}(0)\times [0, 1]} |\Delta u^i|^2\nonumber\\
&\le C\big(\int_{B_{1}(0)}|\nabla u_0^i|^2+\int_{B_{1}(0)\times [0, 1]} |\nabla u^i|^2\big)
\le C(R_i^{\frac2n}+\varepsilon_i^2).
\end{align}

Let $v_0^i:\R^n\to\R^L$ be an $H^1$-extension of $u_0^i$ from $B_{\frac34}(0)$ be $\R^n$ such that
$$v^i_0=u_0^i\ {\rm{on}}\ B_\frac34(0), \ v^i_0=0\  {\rm{in}}\ \R^n\setminus B_1(0), \
\int_{\R^n}|\nabla v_0^i|^2\le C\int_{B_1(0)} |\nabla u_0^i|^2.$$
Define $\widehat{v_0^i}(\cdot, t)=e^{t\Delta} (v_0^i)$.
Then by the local energy inequality we have
\begin{align}\label{LEI50}
&\sup_{0\le t\le 1}\int_{B_{\frac{1}2}(0)}|\nabla \widehat{v_0^i}(t)|^2+\int_{B_{\frac{1}2}(0)\times [0, 1]} |\Delta \widehat{v_0^i}|^2\nonumber\\
&\le C\int_{\R^n}|\nabla {v_0^i}|^2\le C\int_{B_{1}(0)}|\nabla {u_0^i}|^2\le CR_i^{\frac2n}.
\end{align}

Define the blow-up sequence $v_i$ by letting
$$v_i=\frac{u^i-\widehat{v_0^i}}{\varepsilon_i}: B_1(0)\times [0,1]\to \R^L.$$ 
Then $v_i$ solves
\begin{equation} \label{blowup-eqn}
\begin{cases}
\partial_t v_i-\Delta v_i =\varepsilon_i^{-1}A(u^i)(\nabla u^i,\nabla u^i) \ {\rm{in}} \ B_1(0)\times [0,1],\\
v_i\big|_{t=0}=0 \ \ {\rm{on}}\ \ B_\frac34(0),
\end{cases}
\end{equation}
and
\begin{align}\label{h1-bound}
&\int_{B_\frac12(0)\times [0,1]}(|\partial_t v_i|^2+ |\nabla v_i|^2)
\le C\varepsilon_i^{-2}\int_{B_\frac12(0)\times [0,1]}\big(|\nabla u^i|^2+|\nabla \widehat{v_0^i}|^2
+|\Delta u^i|^2+|\Delta \widehat{v_0^i}|^2\big)\nonumber\\
&\le C\varepsilon_i^{-2}\big(R_i^{\frac{2}{n}}+\varepsilon_i^2\big)\le C,
\end{align} 
where we have used \eqref{r-est} in the last step. And it follows from \eqref{small3} and \eqref{r-est} that
\begin{equation}\label{h1-bound1}
\int_{B_\frac12(0)\times [0,1]}|\nabla v_i|^2\le 1+2(1+\varepsilon_i^{-2}R_i^{\frac2n})\le 3+\frac{2C}{i\theta^n}\le 4, 
\end{equation}
and
\begin{equation}\label{h1-bound2}
\theta^{-n}\int_{B_\theta(0)\times [0,\theta^2]} |\nabla v_i|^2> \frac12-\frac{CR_i^{\frac2n}}{\theta^n \varepsilon_i^2}\ge \frac14,
\end{equation}
provided $i$ is sufficiently large.
Furthermore, it follows from \eqref{LEI40} and \eqref{LEI50} that
\begin{equation}\label{LEI60}
\sup_{0\le t\le 1}\int_{B_{\frac{1}2}(0)}|\nabla v_i(t)|^2+\int_{B_{\frac{1}2}(0)\times [0, 1]} |\Delta v_i|^2
\le C\varepsilon_i^{-2}(R_i^{\frac2n}+\varepsilon_i^2)\le C.
\end{equation}

From \eqref{h1-bound}, we may assume that there exists $v_\infty\in H^1(B_\frac12(0)\times [0,1],\R^L)$, with $v(\cdot, 0)=0$ on $B_\frac12(0)$,
such that after passing to a subsequence, $v_i\rightharpoonup v_\infty$ in $H^1(B_\frac12(0)\times [0,1])$. Since
$$\|\varepsilon_i^{-1}A(u^i)(\nabla u^i,\nabla u^i)\|_{L^1(B_\frac12(0)\times [0,1])} \le C\varepsilon_i\to 0,$$
we have that
\begin{equation}\label{blowup-limit}
\partial_t v_\infty-\Delta v_\infty =0\ \ {\rm{in}}\ \ B_\frac12(0)\times [0,1], \ \ v_\infty(\cdot,0)=0 \ \ B_\frac12(0). 
\end{equation}
By the regularity theory of heat equation, we see that
\begin{equation}\label{decay2}
\theta^{-n}\int_{B_\theta(0)\times [0,\theta^2]}|\nabla v_\infty|^2\le C\theta^2\int_{B_\frac12(0)\times [0,1]}|\nabla v_\infty|^2\le C\theta^2.
\end{equation}
On the other hand,  the bound \eqref{LEI60} actually implies that
$v_i\rightarrow v_\infty$ in $H^1(B_\frac12(0)\times [0,1])$. Thus  \eqref{decay2} contradicts to \eqref{h1-bound2}, provided
$\theta\in (0,1)$ is chosen sufficiently small. This completes the proof of \eqref{energy-decay10}.
\qed

\bigskip
With the help of Lemma \ref{energy-decay1} and \eqref{normalized-2-est-u}, we can prove Theorem \ref{initial-Holder} as follows.

\medskip
\noindent{\bf Proof of Theorem \ref{initial-Holder}}:  First, observe that  if we choose $\delta>0$ so small
that $C\delta^2\le \varepsilon_0^2,$
then \eqref{normalized-2-est-u} implies that  \eqref{small1} holds for any $x_0\in\R^n$ and $R_0>0$. 

We can repeatedly apply Lemma \ref{energy-decay1} $k$-times to obtain that for any $x_0\in\R^n$, with $\frac12\le |x_0|\le 2$, 
\begin{equation*}
(\theta_0^kR_0)^{-n}\int_{B_{\theta_0^kR_0}(x_0)\times [0,(\theta_0^kR_0)^2]}|\nabla u|^2\le \frac1{2^k}\max\Big\{
R_0^{-n}\int_{B_{R_0}(x_0)\times [0,R_0^2]}|\nabla u|^2,\ \frac{C_0 R_0^{\frac2n}}{1-\theta_0^{\frac{2}n}}\Big\}.
\end{equation*}
Set $\alpha_0=\min\big\{\frac{\ln 2}{2|\ln \theta_0|}, \frac1n\big\}\in (0,1)$.  Then we obtain that for any  $x_0\in\R^n$, with $\frac12\le |x_0|\le 2$,
and $0<r\le R_0$, it holds that
\begin{equation}\label{energy-decay12}
r^{-n}\int_{B_{r}(x_0)\times [0, r^2]}|\nabla u|^2\le \big(\frac{r}{R_0}\big)^{2\alpha_0} \max\Big\{
R_0^{-n}\int_{B_{R_0}(x_0)\times [0,R_0^2]}|\nabla u|^2,\ \frac{C_0 R_0^{\frac2n}}{1-\theta_0^{\frac{2}n}}\Big\}.
\end{equation}
This, combined with \eqref{LEI20}, further implies that 
for any  $x_0\in\R^n$, with $\frac12\le |x_0|\le 2$,
and $0<r\le R_0$, it holds that
\begin{align}\label{energy-decay13}
&r^{-n}\int_{B_{r}(x_0)\times [0, r^2]}\big(r^2|\partial_t u|^2+|\nabla u|^2\big)\nonumber\\
&\le \big(\frac{r}{R_0}\big)^{2\alpha_0} \max\Big\{
R_0^{-n}\int_{B_{R_0}(x_0)\times [0,R_0^2]}|\nabla u|^2,\ \frac{C_0 R_0^{\frac2n}}{1-\theta_0^{\frac{2}n}}\Big\}.
\end{align}

Now for any $\displaystyle 0<t_0\le \frac{R_0^2}{16}$,   it follows from Proposition 2.5
and \eqref{energy-decay13} 
that\\
i) When $0<r<\sqrt{t_0}$, it follows from the interior gradient estimate by Proposition \ref{interior-smooth} that
\begin{align}\label{interior-energy-decay1}
&r^{-n}\int_{B_{r}(x_0)\times [t_0-r^2, t_0+r^2]}\big(r^2|\partial_t u|^2+|\nabla u|^2\big)\nonumber\\
&\le C\big(\frac{r}{\sqrt{t_0}}\big)^2t_0^{-\frac{n}2}\int_{B_{\sqrt{t_0}}(x_0)\times [0, 2t_0]}|\nabla u|^2\nonumber\\
&\le C\big(\frac{r}{R_0}\big)^{2\alpha_0}\max\Big\{
R_0^{-n}\int_{B_{R_0}(x_0)\times [0,R_0^2]}|\nabla u|^2, \ \frac{C_0 R_0^{\frac2n}}{1-\theta_0^{\frac{2}n}}\Big\}.
\end{align}
ii) When $\displaystyle\sqrt{t_0}\le r\le \frac{R_0}4$, by \eqref{energy-decay13} we also have that
\begin{align}\label{interior-energy-decay2}
&r^{-n}\int_{B_{r}(x_0)\times [0, t_0+r^2]}\big(r^2|\partial_t u|^2+|\nabla u|^2\big)\nonumber\\
&\le 2^{\frac{n}2}(\sqrt{2}r)^{-n}\int_{B_{\sqrt{2} r}(x_0)\times [0, 2r^2]}\big(r^2|\partial_t u|^2+|\nabla u|^2\big)\nonumber\\
&\le C\big(\frac{r}{R_0}\big)^{2\alpha_0}\max\Big\{
R_0^{-n}\int_{B_{R_0}(x_0)\times [0,R_0^2]}|\nabla u|^2, \ \frac{C_0 R_0^{\frac2n}}{1-\theta_0^{\frac{2}n}}\Big\}.
\end{align}

Putting together \eqref{energy-decay13} with \eqref{interior-energy-decay1} and \eqref{interior-energy-decay2}, and applying Morrey's decay Lemma
\cite{Morrey2008}, we conclude that $u\in C^{\alpha_0}\big((B_2\setminus B_\frac12)\times [0, \frac{R_0^2}{16}], N)$ and
\begin{equation}\label{initial-Holder1}
\big[u\big]_{C^{\alpha_0}\big((B_2\setminus B_\frac12)\times [0, \frac{R_0^2}{16}])\big)}
\le C\big(\varepsilon_0, R_0, \|\nabla u_0\|_{L^n(\mathbb S^{n-1})}\big).
\end{equation}

Now we sketch how to improve $\alpha_0$ to $\frac1n$ as follows. To do this, 
let $v_0\in C^{\frac1n}(\R^n, \R^L)$ be such that 
\begin{equation*}
\begin{cases}
v_0=u_0 \ \ {\rm{in}}\ \ B_2\setminus B_\frac12; \ \ v_0=0\ \ {\rm{outside}}\ \ B_3\setminus B_\frac14, \\
\big\|v_0\big\|_{C^0(\R^n)}\le C\big\|u_0\big\|_{C^0(B_2\setminus B_\frac12)}\le C;\\
\ \big[v_0\big]_{C^{\frac1n}(\R^n)}\le C \big[u_0\big]_{C^{\frac1n}(B_2\setminus B_\frac12)}\le C\big\|\nabla_{\mathbb S^{n-1}}u_0\big\|_{L^n(\mathbb S^{n-1})}.
\end{cases}
\end{equation*}
Set $\widehat{v_0}(t)=e^{t\Delta} v_0$. Then
\begin{equation*}
\big\|\widehat{v_0}(t)\big\|_{C^0(\R^n)}\le \big\|v_0\big\|_{C^0(\R^n)}\le C; \ \big[\widehat{v_0}(t)\big]_{C^{\frac1n}(\R^n)}\le 
\big[v_0\big]_{C^{\frac1n}(\R^n)}\le
C, \ \forall t>0,
\end{equation*}
and
\begin{equation}\label{gamma-cont}
R^{-n}\int_{B_R(x_0)\times [0, R^2]} |\nabla \widehat{v_0}|^2
\le C[v_0]_{C^{\frac1n}(\R^n)}^2 R^{\frac{2}{n}} \le C R^{\frac2n}, \ \forall x_0\in \R^n, \ R>0.
\end{equation}
Now for any fixed $x_0\in B_\frac{15}{8}\setminus B_\frac58$ and $0<R<\frac18$, let $w: B_{R}(x_0)\times [0, R^2]\to\R^L$ solve
\begin{equation*}
\begin{cases}
\partial_t w-\Delta w=0\ \ {\rm{in}}\ \ B_{R}(x_0)\times [0, R^2], \\
w(\cdot, 0)=0\ \ {\rm{on}}\ \ B_R(x_0); \ \ w=u-\widehat{v_0} \ \ {\rm{on}}\ \ \partial B_R(x_0)\times [0,R^2].
\end{cases}
\end{equation*}
By the maximum principle and the interior gradient estimate of heat equations, we have that
\begin{equation}\label{w-bound1}
\big\|w\big\|_{L^\infty(B_R(x_0)\times [0, R^2])} \le \big\|u-\widehat{v_0}\big\|_{L^\infty(\partial B_R(x_0)\times [0, R^2])}\le CR^{\alpha_0},
\end{equation}
and
\begin{equation}\label{w-bound2}
\displaystyle r^{-n}\int_{B_r(x_0)\times [0, r^2]} |\nabla w|^2\le C(\frac{r}{R})^2 R^{-n}\int_{B_R(x_0)\times [0, R^2]}|\nabla w|^2,
 \ \forall 0<r\le R.
\end{equation}
Multiplying \eqref{HF1} by $u-\widehat{v_0}-w$ and integrating the equation
over $B_R(x_0)\times [0, R^2]$ yields
\begin{align}\label{w-bound3}
\int_{B_R(x_0)\times [0, R^2]} 
|\nabla(u-\widehat{v_0}-w)|^2&\le C\int_{B_R(x_0)\times [0, R^2]} |\nabla u|^2 |u-\widehat{v_0}-w|\nonumber\\
&\le C\big\|u-\widehat{v_0}-w\big\|_{L^\infty(B_R(x_0)\times [0, R^2])} \int_{B_R(x_0)\times [0, R^2]} |\nabla u|^2\nonumber\\
&\le CR^{\alpha_0} \int_{B_R(x_0)\times [0, R^2]} |\nabla u|^2\nonumber\\
&\le CR^{n+3\alpha_0},
\end{align}
where we have applied \eqref{energy-decay12} in the last step. 

Combining \eqref{w-bound3} with \eqref{gamma-cont} and \eqref{energy-decay12}, we also have
\begin{equation}\label{w-bound4}
R^{-n}\int_{B_R(x_0)\times [0, R^2]}|\nabla w|^2\le CR^{2\alpha_0}.
\end{equation}
Substituting \eqref{w-bound4} into \eqref{w-bound2} and employing \eqref{w-bound3} and \eqref{gamma-cont}, we conclude
that for any $0<r<\frac{R}2$, 
\begin{equation}\label{energy-decay20}
r^{-n}\int_{B_{r}(x_0)\times [0, r^2]}|\nabla u|^2
\le C(\frac{r}{R})^2  R^{2\alpha_0} +Cr^{-n} R^{n+3\alpha_0}.
\end{equation}
Let $\displaystyle\beta_0=\frac{\alpha_0}{n+2}\in (0,1)$.
Then $\displaystyle\alpha_1=\frac{\alpha_0+\beta_0}{1+\beta_0}=\frac{(n+3)\alpha_0}{n+2+\alpha_0}
\in (\alpha_0,1)$. And \eqref{energy-decay30} implies that
\begin{equation}\label{energy-decay30}
\big(R^{1+\beta_0}\big)^{-n}\int_{B_{R^{1+\beta_0}}(x_0)\times [0, R^{2(1+\beta_0)}]}|\nabla u|^2
\le C\big(R^{1+\beta_0}\big)^{2\alpha_1}
\end{equation}
holds for any $x_0\in B_\frac{15}{8}\setminus B_\frac58$ and $0<R<\frac18$. Repeating the same arguments
as in \eqref{energy-decay13}, \eqref{interior-energy-decay1}, and \eqref{interior-energy-decay2}, we can conclude
that $u\in C^{\alpha_1}((B_\frac{15}{8}\setminus B_\frac58)\times [0, \frac{1}{16}])$. It is not hard to see that this process
will take at most finitely many steps until the H\"older exponent $\alpha_0$ of $u$ in $(B_\frac32\setminus B_\frac14)\times [0, \frac{1}{16}]$
reaches $\frac1n$. 
This, combined with the forward self-similarity, implies that
$u\in C^{\frac1n}( \R^n\times [0,\frac{1}{16}]\setminus \{(0,0)\})$. 
This completes the proof of Theorem \ref{initial-Holder}.\qed

\bigskip
\noindent{\bf Acknowledgments:}  
The first author is partially supported by an AMS-Simons travel grant. 
The second author is partially supported by NSF grant 2101224.

\bigskip

\noindent{\bf Data availability}  Data sharing not applicable to this article as no datasets were generated
or analysed during the current study.

\bigskip
\noindent{\bf Conflict of interest} There are no Conflict of interest with third parties.

\end{document}